\documentclass[11pt]{article}
\usepackage{amscd,amssymb,amsmath,verbatim,amsthm}
\usepackage{hyperref}
\usepackage{times}
\usepackage{enumerate}
\usepackage[11pt,curve,matrix,arrow,frame,graph]{xy}
\newcounter{intro}

\newtheorem{theo}[intro]{Theorem}

\newtheorem{thm}{Theorem}[section]

\newtheorem{prop}[thm]{Proposition}

\theoremstyle{remark}

\newtheorem{rems}[thm]{Remarks}

\numberwithin{equation}{section}   

\newcounter{counteroman}

\newcommand{\cref}[1]{Corollary~\ref{#1}}

\newcommand{\pref}[1]{Proposition~\ref{#1}}
\newcommand{\rref}[1]{Remark~\ref{#1}}
\newcommand{\tref}[1]{Theorem~\ref{#1}}

\newcommand{\Poin}{\mathrm{Poin}}
\newcommand{\Sob}{\mathrm{Sob}}

\newcommand{\R}{\mathbb{R}}
\newcommand{\C}{\mathbb{C}}

\newcommand{\bB}{\mathbb{B}}

\newcommand{\N}{\mathbb{N}}

\newcommand{\bS}{\mathbb{S}}

\newcommand{\diver}{\mathrm{div}}
\newcommand{\reg}{\mathrm{reg}}

\newcommand{\cN}{\mathcal{N}}

\newcommand{\bpm}{\begin{pmatrix}}
\newcommand{\epm}{\end{pmatrix}}

\DeclareMathOperator{\spec}{spec}

\DeclareMathOperator{\diam}{diam}

\DeclareMathOperator{\vol}{vol}

\newtheorem{theorem}{Theorem}

\newtheorem{definition}[theorem]{Definition}

\newtheorem{remark}[theorem]{Remark}

\numberwithin{equation}{section}
\numberwithin{theorem}{section}

\renewcommand{\tilde}{\widetilde}
\renewcommand{\bar}{\overline}

\newcommand\dvol{\operatorname{dvol}}

\newcommand\paperintro%
        {%
         }
\newcommand\paperbody%
        {%
         }

\newcommand\cC{\mathcal{C}}
\newcommand\cD{\mathcal{D}}
\newcommand\cE{\mathcal{E}}

\newcommand\cW{\mathcal{W}}

\DeclareMathAlphabet{\mathpzc}{OT1}{pzc}{m}{it}

\newcommand{\del}{\partial}

\newcommand{\calC}{{\mathcal C}}

\newcommand{\calO}{{\mathcal O}}

\newcommand{\calU}{{\mathcal U}}

\newcommand{\calW}{{\mathcal W}}

\newcommand{\eucl}{g_{\rm eucl}}

\begin{document}
\title{H\"older regularity of solutions for \\
Schr\"odinger operators on stratified spaces}
\author{Kazuo Akutagawa \thanks{Supported in part by Grant-in-Aid for Scientific Research (B), JSPS, No. 24340008, 
akutagawa@math.titech.ac.jp}\\ Tokyo Institute of Technology \and 
Gilles Carron \thanks{Supported in part by ANR grant ACG: ANR-10-BLAN 0105, Gilles.Carron@univ-nantes.fr} \\ Universit\'e de Nantes  \and 
Rafe Mazzeo \thanks{Supported in part by NSF DMS 1105050; mazzeo@math.stanford.edu} \\ Stanford University}





\maketitle
\begin{abstract}
We study the regularity properties for solutions of a class of Schr\"odinger equations $(\Delta + V) u = 0$
on a stratified space $M$ endowed with an iterated edge metric.  The focus is on obtaining optimal H\"older 
regularity of these solutions assuming fairly minimal conditions on the underlying metric and potential. 
\end{abstract}

\section{Introduction }
Let $(M,g)$ be a smoothly stratified space with an iterated edge metric, and suppose that $V \in L^p(M; \dvol_g)$.
We prove in this paper that any $W^{1,2}$ solution of the Schr\"odinger equation $(\Delta_g + V)u =0$ 
satisfies a H\"older condition of order $\mu$, where $\mu$ is determined by $p$ and the geometry of $(M,g)$.  
When $g$ and $V$ are polyhomogeneous, i.e., admit asymptotic expansions around each singular stratum in 
powers of the distance function to that stratum (this is the appropriate notion of smoothness in the category of stratified
spaces), then it is known that the solution $u$ is also polyhomogeneous. This is proved using the machinery
of geometric microlocal analysis, see \cite{ACM1}. The exponents which appear in the expansions for $u$ are 
determined by global spectral data on the links of the corresponding strata, and are typically not integers. The appearance 
of a term $r^\mu$ with $\mu \in (0,1)$ in such an expansion shows that from a certain perspective, H\"older regularity is 
the best that could be expected. Our goal here is to show that such H\"older regularity results can be obtained more directly 
and with more classical methods, also allowing metrics which are themselves only of limited regularity.  This is quite 
useful in many situations, for example certain nonlinear problems in geometry, where one may not know the optimal 
regularity of the metric $g$ beforehand. 

We begin by recalling briefly the definition of smoothly stratified spaces; details are deferred until \S \ref{iem} below.  
A topological space $M$ is called smoothly stratified if it decomposes into the union of open manifolds $Y_k$ of varying dimensions 
($\dim Y_k = k$, $k = 0, \ldots, n$) which fit together in a precise manner. We assume that the top-dimensional stratum
$\Omega := Y_n$ is open and dense in $M$, and that $Y_{n-1} = \emptyset$, i.e., $M$ has no codimension $1$ boundary. 
The union of strata $Y_k$, $k < n$, is called the singular set $\Sigma$, sometimes also denoted $M^{\mathrm{sing}}$, 
while $\Omega$ is called the regular set $M^{\mathrm{reg}}$. The crucial property is that each stratum has a tubular neighborhood 
$\calU_k$ which is identified with a bundle of truncated cones over $Y_k$, with fibre $C_R(Z_k)$, where the link $Z_k$ of each 
conic fibre is a stratified space of `lower complexity' and the radial variable of the cone lies in $[0,R)$. To elaborate on this 
local identification, for each $x\in Y_k$, there exists a radius $\delta_x > 0$, a neighborhood $\cW_x$ of $x$ in $M$, and a 
homeomorphism 
\begin{equation}
\varphi_x \colon \bB^k(\delta_x)\times C_{\delta_x}(Z_k)\rightarrow \cW_x,
\label{defvarphi}
\end{equation}
which restricts to a diffeomorphism between $(\bB^k(\delta_x)\times C_{\delta_x}(Z^{\mathrm{reg}}))\setminus 
(\bB^k(\delta_x) \times \{0\})$ and $\cW_x\cap M^{\mathrm{reg}}$.

An iterated edge metric $g$ on $M$ is a smooth (or just H\"older continuous) Riemannian metric on $M^{\reg}$ which
is a perturbation of the model product metric $g_0 = \eucl+dr^2+r^2k_Z$ near each stratum, where $k_Z$ is an iterated 
edge metric on the stratified space $Z$. More specifically, for some $\gamma > 0$, $g$ is locally H\"older of 
order $\gamma$ on $M^{\reg}$ with respect to $g_0$ and satisfies 
\begin{equation}
\label{asymp}
\big|\varphi^*g-g_0\,\big|_{g_0} \le C r^{\gamma}, \quad  \mbox{on}\ \ \bB^k(r)\times C_{r}(Z^{\rm reg}) 
\end{equation}
for all $ r < \delta_x$.  

The simplest nontrivial stratified space is one with simple edges. Such a space has only one nontrivial stratum
$Y_k$, and the link $Z_k$ of the corresponding cone-bundle is a smooth compact manifold of dimension $n-k-1$. 
The best known of these are the spaces with isolated conic singularities, i.e., where $k=0$ here. 

Let $\Delta$ be the Laplace operator on $M^{\reg}$ associated to a given iterated edge metric $g$.  There is an unbounded self-adjoint 
operator $-\Delta$ on $L^2(M, \dvol_g)$ obtained by the Friedrichs extension method and associated to the semi-bounded 
quadratic form $\cC^1_0(M^{\reg}) \ni u \mapsto \int_M |du|_g^2\, \dvol_g$.  When $g$ is only H\"older continuous, it is
necessary to regard $-\Delta$ as the abstract self-adjoint operator associated to this quadratic form, which makes good sense, even
though the coefficients of this differential operator are distributional. We proceed with this understanding, but rarely
mention it again. It is proved in \cite{ACM1, ACM2} that in this setting,
the Riemannian volume form $\dvol_g$ is a doubling measure, and there are Poincar\'e and Sobolev inequalities. 
Consequently, adapting Moser's classical method, we showed that if $V\in L^p$ for some $p>n/2$, then a solution $u$ 
of the equation $(\Delta +V) u=0$ lies in a H\"older space of order $\mu$ for some $\mu \in (0,1)$, see \cite[Theorem 4.8]{ACM1}.  
Our goal in this paper is to understand
the optimal H\"older exponent $\mu$; as we shall show that this optimal exponent has a geometric interpretation.

To state our results, first recall that if $(W,h)$ is a compact stratified space with iterated edge metric, and $\dim W = \ell$,
then it is shown in \cite{ACM1} that $-\Delta_h$ has discrete spectrum. Let $\lambda_1(W)$ denote its first nonzero eigenvalue, 
and also define 
\begin{equation}
\nu_1(W) = \begin{cases}  1 \qquad \mathrm{if}\ \lambda_1(W)\geq \ell, \\
\mbox{the unique value in $(0,1)$ such that } \\
\quad \lambda_1(W)=\nu_1(W)\left(\ell -1 +\nu_1(W)\right) \quad  \mathrm{if}\ \lambda_1(W) < \ell.
\end{cases}
\label{defnu1}
\end{equation}

\begin{theo}\label{A} Let $(M^n,g)$ be a smoothly stratified space with an iterated edge metric. For each $x \in M$,
denote by $Z_x$ the link of the cone bundle over the stratum containing $x$, as in \eqref{defvarphi}, and define 
\begin{equation}
\nu(M) =\inf_{x\in M}\nu_1(Z_x).
\label{defnu}
\end{equation}
Now let $u \in W^{1,2}$ be a solution to $\Delta u+ Vu=0$, where $V\in L^p$. 
\begin{enumerate}[i)]
\item If $V\in L^\infty$ and $\nu=1$, then there is a constant $C>0$ such that for all $x,y\in M$ with 
$d_g(x,y) \leq 1/2$, 
\[
\left|u(x)-u(y)\right|\le C \sqrt{|\log d_g(x,y)|}\, d_g(x,y).
\]
\item If $V\in L^\infty$ and $\nu\in (0,1)$, then $u \in \calC^{0, \nu}(M)$. 
\item if $V\in L^p$ for some $p\in (n/2, \infty)$ and $\nu\in (0,1]$, then $u \in \calC^{0,\mu}(M)$, where
\[
\mu=\min \left\{\nu, 1-\frac n{2p}\right\}.
\]
\end{enumerate}
\end{theo}

As explained above, the novelty of this result is that it requires very little regularity on the metric $g$.  It is known, 
see \cite[section 3]{ACM1}, that when $g$ and $V$  are polyhomogeneous, and the operator $\Delta$ has constant 
indicial roots in some range, then the solution $u$ has a partial polyhomogeneous expansion. This stronger result requires 
quite a lot of machinery to prove, whereas the Theorem above is obtained using more general arguments using only \eqref{asymp}. 

In the course of the proof we shall use a description of neighborhoods in $M$ slightly different than the product
decomposition \eqref{defvarphi}. Namely, it follows easily from \eqref{defvarphi} and \eqref{asymp} that at each 
point $x \in M$ there is a unique tangent cone; this is the Gromov-Hausdorff limit of the family of pointed metric 
spaces $(M,\lambda \mathrm{dist}_g , x)$ as $\lambda\nearrow \infty$.   This limit is unique and is an exact metric 
cone $(C(S_x), dt^2 + t^2 h_x)$ over a compact smoothly stratified space $S_x$, called the tangent sphere at $x$, 
where $h_x$ is an iterated edge metric on $S_x$. Comparing with \eqref{defvarphi}, we see that 
\[
C(S_x)=\R^k\times C(Z_x).
\]
Thus $S_x$ is the $k$-fold spherical suspension of $Z_x$, i.e., the product $[0, \pi/2] \times \bS^{k-1} \times Z_x$ with  metric
\begin{equation}
h_x = d\psi^2 + \sin^2 \psi \, g_{\bS^{k-1}} + \cos^2 \psi \, k_{Z_x}.
\label{sphersusp}
\end{equation}
Note that $S_x$ is  ``as complicated'' of a stratified space as $M$ itself. For example, if $M$ has a simple edge of 
dimension $k$, then $S_x$ has a simple edge of dimension $k-1$ (in particular, if $M$ has an isolate conic 
singularity at $x$, then $S_x = Z_x$ is a smooth compact manifold). 

The reason we bring this up now is that much of the analysis below is done on cones $C(S)$, either 
with respect to an exact conic metric $g_0$ or one which is a small perturbation of it. The result of this analysis
is that a solution $u$ of $(\Delta + V) u  = 0$ on $C(S)$ lies in the H\"older class of order $\mu$ where
$\mu$ is determined the H\"older exponent $p$ for $V$ and the constant $\nu(S)$ from \eqref{defnu}. 
To obtain the result above, we must then show that 
\begin{equation}
\mbox{if}\ \ C(S) = \R^k \times C(Z),\ \ \mbox{then} \ \ \nu(S) = \nu(Z). 
\label{nuSnuZ}
\end{equation}
This is proved in \S \ref{spectral} below.

The main step in proving Theorem \ref{A} is to show that under the hypothesis of case iii), $u$ satisfies the Morrey condition
\[
\frac{1}{\vol B(x,r)}\int_{B(x,r)} |du|^2\le Cr^{2-2\mu}
\]
for all $x\in M$ and $r\in (0,1)$.  It is well known that for Dirichlet spaces which are measure doubling and have a 
Poincar\'e inequality, such an estimate yields the H\"older continuity of $u$.  We recall the proof of this in the appendix. 

One difficulty in the analysis is that the comparison between the geometry of $M$ near a point $x$ and of the 
tangent cone $C(S_x)$ can only be made below a certain length scale $\delta_x$. In the next section, we describe some 
facts from the geometry of balls which allow us to circumvent this difficulty. In \S \ref{DN}, we develop some familiar 
analytical tools on stratified spaces, namely the Green formula and the Dirichlet-to-Neumann operator, which are used
in the later analysis. This is followed by a monotonicity formula for the quadratic form associated to $\Delta+V$. 
Theorem \ref{A} is proved in \S \ref{proofs}.

\section{On the geometry of stratified space}
We recall some further aspects of the definition of smoothly stratified spaces, all taken from \cite[\S 2]{ALMP}, and then state some 
facts about the structure of balls and tangent cones for these spaces. We refer to \cite{ALMP} for further details.

\subsection{Stratifications and iterated edge metrics}
\label{iem}
Let $M$ be a smoothly stratified space. As described in the introduction, this means that $M = \sqcup_{j \leq n} Y_j$,
where $Y_j$ is a (typically open) smooth manifold of dimension $j$.  We assume that $M$ is compact and
$Y_{n-1} = \emptyset$. Any $x \in Y_j$ has a neighborhood homeomorphic to
$\bB^{j}(\eta) \times C_\eta (Z)$, where $Z$ is a stratified space of dimension $n-j-1$, $C_\eta(Z)$ is the metric cone 
over $Z$ truncated at radius $\eta$ and $\bB^{j}(\eta)\subset \R^j$ is a Euclidean ball of radius $\eta$.

The \textit{depth } of a stratum $Y$ is the largest integer $k$ such that there is a chain of strata $Y_{j_1}, \ldots, Y_{j_k}$ with
$Y_{j_{i-1}} \subset \overline{Y_{j_i}}$ and $Y_{j_1} = Y$.  A stratum of maximal depth is necessarily a closed manifold.

The stratified space $M$ can be covered by a finite number of open set $\mathcal{W}_\alpha$, each homeomorphic to 
$\mathcal{U}_\alpha\times C_{\delta_\alpha}(Z_\alpha)$, where, for some $\gamma\in (0,1]$,  
\begin{itemize}
\item $\mathcal{U}_\alpha$ is an open set in $\R^{d_\alpha}$ endowed with a smooth Riemannian metric $h_\alpha$;
\item $Z_\alpha$ is a compact stratified space of dimension $n-d_\alpha-1$ endowed with a uniform $\gamma$-H\"older family 
of iterated edge metric 
\[
\left\{k_\alpha(y),y\in \mathcal{U}_\alpha\,\right\}\,\, ;
\]
\item $C_{\delta_\alpha}(Z_\alpha)$ is the cone over $Z_\alpha$ truncated at radius $\delta_\alpha$; this cone is also a stratified space.
\item $|\varphi^* g - (h_\alpha + dr^2 + r^2 k_\alpha)| \leq C r^\gamma$.
\end{itemize}
We assume that the family of quadratic forms 
\[
\left\{h_\alpha(y)=\sum_{i,j} h_{\alpha,i,j}(y) dy_idy_j\,\,,\,\, y\in \mathcal{U}_\alpha\right\}
\]
is uniformly $\gamma$-H\"older  and precompact, so in particular there are positive constants $c, C$ such that for all $y, y_0\in  \mathcal{U}_\alpha$,
\[
c\, h_\alpha(y_0)\le h_\alpha(y)\le C\, h_\alpha(y_0).
\]

\subsection{The geometry of geodesic balls}
We now describe the geometry of balls $B(m,\tau)$ in $(M,g)$. The main conclusion is that these balls look like
truncated cones $C_\tau(S)$ with a uniformly controlled error. 
 
Choose $\eta>0$ sufficiently small so that any geodesic ball of radius $\eta$ lies in one of the open sets $\mathcal{W}_\alpha$.
Let $\delta\in (0,1)$ be a parameter whose value will be specified below. We study geodesics balls $B(m,\tau)$,
where $\tau\in (0,\delta\eta/4)$. For each such ball, choose an open set $\mathcal{W}_\alpha$ which contains it, and write
$\varphi_\alpha\colon  \mathcal{W}_\alpha\rightarrow \mathcal{U}_\alpha\times C_{\delta_\alpha}(Z_\alpha)$ for the homeomorphism 
$\mathcal{W}_\alpha \to\mathcal{U}_\alpha\times C_{\delta_\alpha}(Z_\alpha)$.  Thus $m \in \calW_\alpha$ has coordinates
$\varphi_\alpha(m)=(y,\rho,z)$.
 \begin{description}
\item[Case 1: $\rho\le \tau/\delta$:]  Setting $\underline{m}=\varphi_\alpha^{-1}\left(y,0,z\right)$, then by the triangle inequality
\[
B(m,2\tau)\subset B\left(\underline{m}, \left(1+1/\delta \right)2\tau\right)\,.
\]
We wish to compare the metric $g$ on this latter ball to the model product metric
\[
g_0=h_\alpha(y)+dr^2+r^2 k_\alpha(y)\, .
\]
Clearly, if $\varepsilon < \eta/2$, then
\[
\left|g-g_0\right|\le C \varepsilon^\gamma \,\, \mathrm{on}\,\, B\left(\underline{m}, \varepsilon \right)\, .
\]

There is a constant $\kappa$ such that the $\bB_0(\kappa r)$, which is the same as the cone $C_{\kappa r}(S_{\underline{m}})$
satisfies 
\[
B(m,\tau)\subset \bB_0(\kappa\tau)\subset B(m,2\kappa\tau).
\] 
Furthermore, on $B(m,\tau)$, we have
\[
\left|g-g_0\right|\le C\tau^\gamma.
\]

\item[Case 2: $\rho\ge \tau/\delta$: ] On $B(m,2\tau)$ we have
\[
g=h_\alpha(y)+ dr^2+\rho^2k_\alpha(y)+ \calO\left(\frac{ \tau^\gamma}{ \rho^\gamma} \right).
\] 
Furthermore, if $\delta$ is small enough, 
\[
\varphi_\alpha(B(m,2\tau))\subset \bB(x,3\tau)\times (\rho-3\tau,\rho+3\tau)\times B^{Z_\alpha}\left(z,3 \tau/\rho\right),
\]
where $B^{Z_\alpha}\left(z,3 \tau/\rho\right)$ is the ball of radius $3\tau/\rho$ in $\left(Z_\alpha,k_\alpha(y)\right)$.
\end{description}

Using the relationships and estimates in these two cases, we can then prove the following, via an induction on the depth of the
stratified space. 
\begin{prop}\label{geoBoule}
There are positive constants $\Lambda, \delta_0,\kappa$, with $\Lambda \delta_0<1$, such that for any $\delta\in (0,\delta_0)$
and $m\in M$, if $B(m,\delta)\subset \mathcal{W}_\alpha$, then there is a sequence of numbers
$\rho_1 > \rho_2 > \dots > \rho_{d_\alpha} > \rho_{d_{\alpha+1}}=0$ so that if we set $\tau_j=\delta \,\prod_{i=1}^j \rho_j$ and 
choose any $\tau \in [\rho_j, \rho_{j-1})$, then there is an open set $\Omega_{m,j,\alpha}$ homeomorphic to a cone
$C_{\kappa \tau}(S_{m,j,\alpha})$ over a connected stratified space $S_{m,j,\alpha}$ such that
\[
B(m,\tau)\subset \Omega_{m,j,\alpha}\subset B(m,2\kappa \tau).
\]
Moreover there is an iterated edge metric $h_{m,j,\alpha}$ on $S_{m,j,\alpha}$ so that on $\Omega$,
\[
\left|g-\left( dt^2+t^2 h_{m,j,\alpha}\right)\right|\le \Lambda\,\left(\frac{\tau}{\rho_1\rho_2\dots\rho_{j-1}}\right)^\gamma\,\, .
\]
The set of metric spaces $(S_{m,j,\alpha},h_{m,j,\alpha})$, where $m,j$ and $\alpha$ vary, for a fixed $(M,g)$, is precompact in 
the biLipschitz topology on the space of all compact metric spaces. In particular there is a finite set of compact metric spaces 
$(Y_j,d_j)$, $j = 1, \ldots, N$, and a constant $K>1$ so that each $(S_{m,j,\alpha},h_{m,j,\alpha})$ is  $K$-biLipschitz to 
at least one of the $(Y_\beta,d_\beta)$. 
\end{prop}


\section{Some analytical tools} \label{DN}

\subsection{The Poincar\'e and Sobolev inequalities}
We first recall briefly the proof that any compact stratified space with iterated edge metric satisfies a scale-invariant Poincar\'e 
inequality, and hence also a Sobolev inequality, and hence the Laplace operator has discrete spectrum. We prove this first 
under a topological condition, but then explain in a remark how this condition may be removed.  

\begin{prop}
Let $(M,g)$ be a compact stratified space with an iterated edge metric. Assume that for each $x\in M$, the tangent 
sphere $S_x$ is connected.  Then there are constants $a>1$ and $C,\eta>0$ such that if $B$ is any ball of radius $r(B) < \eta$,
then for every $f\in W^{1,2}(aB)$, there is a scale-invariant Poincar\'e inequality
\[ 
\int_B \left| f-f_B\right|^2 \, \dvol_g \le C_{\Poin} r(B)^2 \int_{aB} \left|d f\right|^2 \, \dvol_g. 
\]
\end{prop}
\begin{remark}\label{cut} 
If the connectedness condition fails for the tangent spheres $S_x$ along certain of the strata, then we can define a new stratified
space $\widetilde{M}$ for which this condition does hold as follows: cut $M$ along each stratum where the corresponding link
is disconnected. The connectedness condition holds for this new space, and the Poincar\'e inequality on $\widetilde{M}$ 
implies one on $M$ as well.
\end{remark}
\begin{remark}  It is known, cf.\ \cite{LSC1}, \cite{LSC2}, \cite[Theorem 5.1]{HK}, that if $(M,g)$ is a space with a scale-invariant 
Poincar\'e inequality, and is such that the measure $dV_g$ is Ahlfors $n$-regular, i.e., $cr^n\le \vol B(x,r)\le C r^n$ for all
$x \in M$ and all $0 <  r < \frac12 \diam_g M$, then there is a Sobolev inequality
\begin{equation}
C_{\Sob} \|\psi\|^2_{L^{\frac{2n}{n-2}}}\le \|d\psi\|^2_{L^2}+\|\psi\|^2_{L^2} \, ,
\label{Sobineq}
\end{equation}
for every $\psi\in W^{1,2}(M)$.   This Sobolev inequality implies, in turn, that the spectrum of the Friedrichs realization of
the Laplace operator $-\Delta_g$ is discrete, i.e.\ there exist $\lambda_j \nearrow \infty$ and $\varphi_j$,  such that
$-\Delta_g \varphi_j = \lambda_j \varphi_j$ and so that the closed linear span of the $\varphi_j$ equals $L^2(M)$. 
\end{remark}

\medskip

The proof of the Proposition is inductive.  We assume that the result has been proved for all compact stratified spaces
(with iterated edge metrics) of depth less than $d$ and then prove that it holds for spaces of depth $d$.  By an obvious
localization argument, it suffices to show that if there is a scale-invariant Poincar\'e inequality on a
connected stratified space $S$, then there is also one on the truncated cone $C_R := C_R(S)$, 
\[ 
\|f-f_{C_R}\|^2_{L^2(C_R)}\le C_{\Poin} R^2  \|df \|^2_{L^2(C_R)}
\]
for all $f\in W^{1,2}(C_R)$. Once we have established this inequality on $C_R$, it then follows that it holds on any compact
stratified space $(M,g)$ of depth $d$, and we then also obtain the Sobolev inequality and discreteness of the spectrum 
of $-\Delta_g$ on all such spaces. This completes the next step of the induction. 

Thus it remains to prove that the scale-invariant Poincar\'e inequality holds on $C_R$, which we do by noting that it suffices
to take $C_{\Poin} = \max \{ A^{-1}, B^{-1}\}$, where 
\begin{itemize}
\item $A$ is the first nonzero eigenvalue of the operator $-\frac{d^2}{dr^2}-\frac{n-1}{r} \frac{d}{dr}$ on $L^2([0,1])$ with
Neumann conditions at $r=1$.  Equivalently, $\sqrt{A}$ is the first positive zero of $(r^{1-\frac{n}{2}} J_{\frac{n}{2}}(r))'$ where $J_\zeta$ 
is the Bessel function of order $\zeta$, and  
\item $B$ is the lowest eigenvalue of the  $-\frac{d^2}{dr^2}-\frac{n-1}{r} \frac{d}{dr}+\frac{\lambda_1}{r^2}$ on $L^2([0,1])$ again
with Neumann conditions at $r=1$, i.e., $\sqrt{B}$ is the first positive zero of $(r^{1-\frac{n}{2}} J_{\nu}(r))'$,  where 
\[ 
\nu=\sqrt{\lambda_1+\left(\frac{n-2}{2}\right)^2},
\]
where $\lambda_1$ is the first nonzero eigenvalue of $-\Delta_h$ (recalling that since $S$ is connected, $\lambda_0=0<\lambda_1$).
\end{itemize}

\subsection{Restriction to the link}  Let $(S,h)$ be a compact, connected smoothly stratified space of dimension $n-1$ 
with iterated edge metric, and consider the cone $(C(S), g_0=dr^2+r^2h)$. 
The ball $\bB_0(\rho)$ centered at $0$ is simply the truncated cone $C_\rho(S)$.

Write the eigenvalues of $-\Delta_h$ as $\lambda_j =\nu_j(n-2+\nu_j)$, which gives the nondecreasing sequence $\{\nu_j\}$.
Since $S$ is connected, $\nu_0=0<\nu_1\le \nu_2\le \dots$. 

There is a restriction map $R$ from the standard Sobolev space $\displaystyle W^{1,2}(\bB_0(\rho))=\left\{\varphi\in L^2(\bB_0(\rho)), 
d\varphi\in L^2 \right\}\,$ to $\del \bB_0((\rho))$. It is easy to check, using the eigenfunction expansion on $S$, that
\[
R\colon W^{1,2}(\bB_0(\rho))\rightarrow H^{1/2}(\partial\bB_0(\rho))),
\]
where 
\[
H^{1/2}(\partial\bB_0(\rho)))=\left\{\sum_j c_j \varphi_j\,,\, \sum_j \nu_j |c_j|^2<\infty\right\}\,\,.
\]

\subsection{Green's formula}
Let $(M,g)$ be a stratified space with iterated edge metric, and suppose that $X$ is a vector field defined
on $M^{\mathrm{reg}}$. The function $\diver_g X$ is defined in the usual way on this regular part, and we say that 
$\diver_g X \in L^1$ if there is a function $\varphi\in L^1(M)$ such that 
\[
\int_M X u \, \dvol_g=\int_M \varphi u \, \dvol_g.
\]
for all $u\in \cC^1_0(M^{\reg})$. Note that $\diver_g X$ only depends on $\dvol_g$, hence if $\tilde g$ is another Riemannian 
metric such that $\dvol_{\tilde{g}} = J \dvol_g$ for some Lipschitz function $J$, then
\[
\diver_{\tilde{g}} X = \diver_g X + X \cdot \nabla \log J.
\]

Our goal in this subjection is to establish Green's formula on $M$ under low regularity assumptions on $X$. We do this first
when $M$ has no codimension $1$ boundary, and then when $M$ is a truncated cone. 

\begin{prop} Let $X$ be an $L^2$ vector field on $M$ such that $\diver_g X \in L^1$. Then 
\[
\int_M \diver_g X \, \dvol_g = 0. 
\]
\end{prop}
\proof  It is standard that this formula holds if $X$ has compact support in $M^{\reg}$.  Thus if $u$ is Lipschitz
with compact support in $M^{\reg}$, then 
\begin{equation}
0=\int_M\mathrm{div}_g(u X) \, \dvol_g=\int_M  Xu \, \dvol_g+\int_M u \diver_g  X\,  \dvol_g.
\label{cutoffgreen}
\end{equation}
Because the volume of the tubular neighborhood of radius $R$ around $M^{\mathrm{sing}}$ is $\calO(R^2)$, we can choose a sequence $\psi_\ell \in \calC^\infty_0(M^{\reg})$ such that $0\leq \psi _\ell\le 1$, $\lim_{\ell} \psi_\ell(x)=1$
for a.e.\ $x$, and $\lim_{\ell}\|d\psi_\ell\|_{L^2}=0$. Using the inequality 
\[
\left|\int_M  X\psi_\ell \,\dvol_g\right|\le \|X\|_{L^2}\|d\psi_\ell\|_{L^2},
\]
the result follows immediately from \eqref{cutoffgreen}. 
\endproof

Let us now turn to the analog of this result on $\bB_0(\rho)$, the truncated metric cone $C_{[0,\rho)}(S)$, with exact conic metric $g_0$,
and with $n \geq 2$.  Since $\dvol_0=r^{n-1}dr \dvol_h$, the volume form on $\del \bB_0(\rho)$ is $d\sigma_0=\rho^{n-1} \dvol_h$. 
\begin{prop}\label{GreenBy} Let $X$ be an $L^2$ vector field on $\bB_0(\rho)$ with $\mathrm{div}_gX\in L^2$; then 
$X_r=\left\langle X,\frac{\partial}{\partial r}\right\rangle=dr(X) \in H^{-1/2}(\partial \bB_0(\rho))$ and for all
$u\in W^{1,2}(\bB_0(\rho))$ we have
\begin{equation}
\int_{\bB_0(\rho)}u\diver_{g_0} X \, \dvol_0+ \int_{\bB_0(\rho)}Xu\, \dvol_0=\int_{\partial\bB_0(\rho)}X_ru\, d\sigma_0. 
\end{equation}
\end{prop}
\proof
Let  $Y=uX$ where $u$ is Lipschitz.  For $f\in \cC^1_0(0,\rho)$, set $v(r,\theta) = f(r)$, $\theta \in S$. By the preceding proposition,
\[
0=\int_{\bB_0(\rho)}\diver_{g_0}(vY)\, \dvol_0=\int_{\bB_0(\rho)}v\diver_{g_0}(Y)\, \dvol_0+\int_{\bB_0(\rho)}Yv\, \dvol_0. 
\]
However, $Y v=f'(r)Y_r$, where $Y_r=\left\langle Y,\frac{\partial}{\partial r}\right\rangle=uX_r.$ The function $K(r):=
\int_{\partial \bB_0(r)}Y_r \, d\sigma_0$ is thus in $L^1$, so if we write $I(r)=\int_{ \bB_0(r)}\diver_0 Y \, \dvol_0$, then
\[
\int_0^\rho f'(r) I(r)\, dr = - \int_0^\rho f(r)I'(r)\, dr=\int_0^\rho f'(r)K(r)\, dr.
\]
Since this holds for every $f\in \cC^1_0(0,\rho)$, the function $K$ is equal almost everywhere to a  continuous function and 
$K(r) = I(r) + c$ for some constant $c$ and all $r\in (0,\rho)$. Since $X$ and $\diver_0X \in L^2$, we obtain
\[
|I(r)|\le  o\left(r^{n/2}\right)\,\,\mathrm{and}\,\, \left|\int_0^r K(r)dr\right|\le  o\left(r^{n/2}\right),
\]
and since $n\ge 2$ we see that $c=0$.  

Now, for $\rho\ge t>s>0$,
\begin{multline*}
\int_{\bB_0(t)\setminus \bB_0(s)}u\diver_{g_0}X\, \dvol_0+\int_{\bB_0(t)\setminus \bB_0(s)} Xu\, \dvol_0 \\
=\int_{\partial\bB_0(t)}X_ru \, d\sigma_0-\int_{\partial\bB_0(s)}X_ru\, d\sigma_0.
\end{multline*}
which gives
\[
\left|\int_{\partial\bB_0(t)}X_ru \, d\sigma_0-\int_{\partial\bB_0(s)}X_ru \, d\sigma_0\right| \le 
\varepsilon(t,s) \, \|u\|_{W^{1,2}}, 
\]
where $\varepsilon(t,s) = ||X||^2_{L^2} + ||\diver_0 X||^2$, these norms being  taken over the annular
region $\bB(0,t)\setminus \bB(0,s)$. 

This holds for all Lipschitz functions $u$, and hence, by a density argument, also when $u \in W^{1,2}$ and for almost every 
$\rho\ge t>s>0$.  Let $f_t(\theta)=X_r(t,\theta)$. If $\varphi\in H^{1/2}(S)$, then let $u$ be the harmonic function on 
$\bB_0(t)\setminus \bB_0(s)$ such that $u=\varphi/s^{n-1}$ on $\partial  \bB_0(s)$ and $u=\varphi/t^{n-1}$ 
on $\partial  \bB_0(t)$. This gives
\[
\left|\int_S (f_t-f_s) \varphi d\sigma\right|\le  \epsilon(t,s)\|\varphi\|_{H^{1/2}}.
\]
In other words, the function $t\mapsto f_t$ is continuous 
as a map $(0,\rho] \longrightarrow H^{-1/2}(S)$.  The asserted formula follows easily from this.
\endproof

Comparing with \cite[\S 5]{BBC1}, we obtain 
\begin{prop} Let $g$ be an iterated edge metric on $\bB_0(\rho)$ such that $\dvol_g=J\dvol_0$, where $J$ is
Lipschitz and $J\ge \epsilon>0$.  
Suppose that $u\in W^{1,2}(\bB_0(\rho))$ and $\Delta_gu\in L^2$ i.e., there exists a constant $C > 0$ such that
\[
\left| \int_{\bB_0(\rho))}\langle du,d\varphi\rangle_g \, \dvol_g\right|\le C \|\varphi\|_{L^2} 
\]
for every $\varphi \in  W_0^{1,2}(\bB_0(\rho))$. Then, letting $n_g$ denote the outward unit normal with respect to $g$, 
$n_g u \in H^{-1/2}(S)$,  and if  $v\in W^{1,2}(\bB_0(\rho))$, then 
\begin{equation}
\int_{\bB_0(\rho)}v\Delta_gu \, \dvol_g+\int_{\bB_0(\rho)}\langle dv,du\rangle_g \, \dvol_g=
\int_{\partial\bB_0(\rho)}v\,n_g u \,  d\sigma_g. 
\end{equation}
\end{prop}\proof
This follows from Proposition \ref{GreenBy} with $X=J  \nabla^g u$. Indeed, $\Delta_gu \, \dvol_g=\mathrm{div}_0(X) \, \dvol_0$,
$\langle dv,du\rangle_g\, \dvol_g= X v\, \dvol_0$, and $n_g u= \frac{\partial u}{\partial r}\frac{1}{|dr|_g}$
and $\frac{dr}{|dr|_g}d\sigma_g=Jdrd\sigma_0$.
\endproof

\subsection{The Dirichlet to Neumann operator}
We now develop properties of the Dirichlet to Neumann operator on the truncated cone $C_\rho(S)$, first with respect to
an exact conic metric and then with respect to a more general iterated edge metric on this space.
\subsubsection{The model case } 
Any $v\in H^{1/2}(\partial\bB_0(\rho)))$ has a unique harmonic extension $\cE_{0,\rho}(v)\in W^{1,2}(\bB_0(\rho))$. 
On eigenfunctions $\varphi_j$ on $S$, there is an explicit formula
\begin{equation}\label{harmonic}
\cE_{0,\rho}(\varphi_j)(r,\theta)=\left(\frac{r}{\rho}\right)^{\nu_j} \varphi_j(\theta)\,\,.
\end{equation}
More generally, $\cE_{0,\rho}(v)$ minimizes the Dirichlet energy 
$$
\int_{\bB_0(\rho)} |d\cE_{0,\rho}(v)|^2_0\dvol_0 \le \int_{\bB_0(\rho)} |d u|^2_0\dvol_0\,\,.
$$
amongst all functions $u\in W^{1,2}(\bB_0(\rho))$
for which the restriction $R(u)$ to the boundary equals $v$, 

\paragraph{Definition} The Dirichlet to Neumann operator $\cN_{0,\rho}$ is the bounded operator :
\begin{equation*}\begin{split}
\cN_{0,\rho}\colon W^{1,2}(S)&\rightarrow L^2(S)\\
v&\mapsto \left.\frac{d}{dr}\right|_{r=\rho}\cE_{0,\rho}(v)(r, \cdot),
\end{split}
\end{equation*}
so in particular 
$$
\cN_{0,\rho}\varphi_j=\frac{\nu_j}{\rho}\varphi_j\,.
$$
The operator $\cN_{0,\rho}$ is selfadjoint with compact resolvent. From the variational characterization of the harmonic extension, 
there is a min-max formula for  its eigenvalues:  
\begin{equation}\label{minmax0}
\frac{\nu_j}{\rho}=\max_{\substack{\scriptscriptstyle V\subset W^{1,2}( \bB_0(\rho))\\ \scriptscriptstyle \dim V=j}}\,\,\, \inf_{u\in V^\perp\setminus\{0\}} \,\,\frac{\int_{\bB_0(\rho)} |du|^2_0\dvol_0}{\int_{\partial\bB_0(\rho)} |u|^2_0d\sigma_0}.
\end{equation}

\subsubsection{The general case} 
Let $g$ be another iterated edge metric on $\bB_0(\rho)$ satisfying
$$
|g-g_0|\le \Lambda \rho^\gamma \ll 1.
$$
Suppose that $V\in L^p(\bB_0(\rho))$ for some $p>n/2$, with the bound
$$
\int_{\bB_0(\rho)} |V|^p\dvol_0\le A^p\,\,.
$$
We shall study properties of the operator $\Delta_g + V$. 

The spaces $L^2(\bB_0(\rho))$ and  $W^{1,2}(\bB_0(\rho))$ are the same relative to either of the two metrics $g_0$ and $g$,
and similarly for  $W^{1,2}(\partial \bB_0(\rho))$. 
We write $W_0^{1,2}(\bB_0(\rho))$ for the set of functions in $u \in W^{1,2}$ such that $R(u) = 0$. Recall too that the space of Lipschitz 
function with compact support in $(0,\rho)\times S^{\reg}=\left(\bB_{0}(\rho)\right)^{\reg}$ is dense in $W_0^{1,2}(\bB_0(\rho))$. 
As in \eqref{Sobineq}, see also \cite{ACM1}, there is a Sobolev inequality for both $g_0$ and $g$, i.e., there exists $C_{\Sob} > 0$ so that
\[
C_{\Sob} \|\psi\|^2_{L^{\frac{2n}{n-2}}}\le \|d\psi\|^2_{L^2} \quad  \forall\psi\in W_0^{1,2}(\bB_0(\rho))
\]
relative to either metric. For $\rho$ sufficiently small, the quadratic form 
$$
\psi\mapsto Q_{g,V,\rho}(\psi) ;= \int_{\bB_0(\rho))}|d\psi|_g^2 \dvol_g-\int_{\bB_0(\rho))}V\psi^2 \, \dvol_g
$$
is coercive in $W_0^{1,2}(\bB_0(\rho))$. Indeed, applying the H\"older inequality twice gives
\begin{equation*}
\begin{split}
\int_{\bB_0(\rho)}V\psi^2 \dvol_g&\le \|V\|^2_{L^{\frac{n}{2}}} \|\psi\|^2_{L^{\frac{2n}{n-2}}}\\
&\le \frac{1}{\mu} A\left(\frac{\vol_h(S)}{n}\rho^n\right)^{\frac2n -\frac1p} \|d\psi\|^2_{L^2},
\end{split}
\end{equation*}
which implies that $Q_{g,V,\rho}(\psi) \geq c Q_{0}(\psi)$ (where $Q_0$ is the quadratic form when $g=g_0$ and $V=0$) provided 
$$
C_{\Sob}^{-1} A\vol_h(S)^{\frac2n -\frac1p} \rho^{2 -\frac n p} <  1.
$$
Assuming this condition, then for each $v\in H^{1/2}(\partial\bB_0(\rho)))$ the functional 
$$
R^{-1}(v) \ni \psi \mapsto  \int_{\bB_0(\rho))}\left(\, |d\psi|_g^2 -V\psi^2\,\right) \, \dvol_g
$$
reaches its infimum at a unique function 
$$
\cE_{V,\rho}(v)\in  W^{1,2}(\bB_0(\rho)).
$$
The Euler-Lagrange condition implies that $\cE_{V,\rho}(v)$ satisfies the equations:
\begin{equation}\label{DiriPro}
\left\{\begin{array}{l}
\left(\Delta_g+V\right)\cE_{V,\rho}(v)=0\\
\left.\cE_{V,\rho}(v)\right|_{\partial \bB_0(\rho)}=v. 
\end{array}\right.  \end{equation}

According to the discussion in \S 3.3, there is a Green formula for functions in the domain
$$
\cD(\Delta_g)=\{\psi\in W^{1,2}(\bB_0(\rho)),\Delta_g\psi\in L^2 \}.
$$
Decompose the $g$ unit normal to $\partial \bB_0(\rho)$ as $\vec n_g=\alpha\frac{\partial}{\partial r}+\vec\beta$, where 
$\vec\beta\perp_{g_0}\frac{\partial}{\partial r}$. Clearly,
\begin{equation}\label{normal}
|\alpha-1|+\left|\vec\beta\right|_{g_0}\le C\Lambda \rho^{\gamma}.
\end{equation}
If $\psi\in\cD(\Delta_g)$, then its normal derivative at the boundary, which we denote by $\frac{\partial \psi}{\partial \vec n_g}$, lies
in $H^{-1/2}(\partial \bB_0(\rho))$ and for any $\varphi\in W^{1,2}(\bB_0(\rho))$ 
\begin{equation}\label{Green}
\int_{\bB_0(\rho))}\Delta_g \psi\,\varphi \dvol_g+ \int_{\bB_0(\rho))}\langle d\psi,d\varphi\rangle_g \,\dvol_g=\int_{\partial\bB_0(\rho))}
\frac{\partial \psi}{\partial \vec n_g} \,\varphi \, d\sigma_g.
\end{equation}
Hence for $\psi,\varphi \in\cD(\Delta_g)$, 
\begin{equation}\label{Green2}
\int_{\bB_0(\rho))}\left(   \psi\,\Delta_g\varphi-\Delta_g \psi\,\varphi \,\right)\dvol_g=
\int_{\partial\bB_0(\rho))}\left(\psi \frac{\partial \varphi}{\partial \vec n_g}- \,\frac{\partial \psi}{\partial \vec n_g} \,\varphi \right)\,d\sigma_g.
\end{equation}

We can now define the Dirichlet to Neumann operator associated to the quadratic form
$$
H^{1/2}(\partial\bB_0(\rho)) \ni v \mapsto \int_{\bB_0(\rho)}\left(\, \left|d\cE_{V,\rho}u\right|^2_g-V\left|\cE_{V,\rho}u\right|^2\,\right)\, \dvol_g.
$$

$$\cN_{V,\rho}v:=\frac{\partial }{\partial \vec n_g}\cE_{V,\rho}(v).$$
The operator $\cN_{g, V,\rho}$ is self-adjoint. We indicate below that it has compact resolvent. It is then not hard to see that its 
spectrum has a min-max interpretation:
\begin{equation}\label{minmaxV}
\mu_j=\max_{\substack{\scriptscriptstyle V\subset W^{1,2}(\bB_0(\rho))\\ \scriptscriptstyle \dim V=j}}
\,\,\, \inf_{u\in V^\perp\setminus\{0\}}
 \,\,\frac{\int_{\bB_0(\rho)}\left[\, |du|^2_g-Vu^2\,\right]\dvol_g}{\int_{\partial\bB_0(\rho)} |u|^2 d\sigma_g}.
\end{equation}

\subsection{Comparison of the spectra}
Our next goal is to compare the spectra of the operators $\cN_{0,\rho}$ and  $\cN_{V,\rho}$.
 
The first step involves finding an $L^p$ estimate for the harmonic extension operator $\cE_{0,\rho}$. If 
\[
v=\sum_j c_j \varphi_j \in H^{1/2}(\partial  \bB_0(\rho))
\]
so that
\[
\cE_{0,\rho}v(r,\theta)=\sum_j\left(\frac{r}{\rho}\right)^{\nu_j} c_j \varphi_j(\theta)=e^{-tL}v(\theta),
\]
where $r=e^{-t}\rho$ and $L=\sqrt{-\Delta_h+\left(\frac{n-2}{2}\right)^2}-\frac{n-2}{2}$.
Assume first that $c_0=0$, i.e., $\int_{\partial  \bB_0(\rho)} v \, d\sigma_0=0$. The Sobolev inequality 
on the product $\left( (1/2,1)\times S, (dr)^2+h\right)$ implies an estimate on the heat kernel of $\Delta_h$. 
The subordination identity 
$$
e^{-tL}= \int_0^\infty \frac{t}{2\sqrt{\pi}} \, e^{\frac{n-2}{2}t-\frac{t^2}{4\tau}-\tau\left(-\Delta_h+\left(\frac{n-2}{2}\right)^2\right)}
\, \frac{d\tau}{\tau^{3/2}}
$$
then shows that if $q\ge 2$ then using $c_0 = 0$,  
$$
\left\| e^{-tL}Lv\right\|_{L^q}\le 
\frac{C}{t^{(n-1)\left(\frac12-\frac1q\right)}}\,\left\| \sqrt{L}v\right\|_{L^2}. 
$$
However, if $q<2(n-1)$, there is also a Sobolev inequality 
$$
\left\| e^{-tL}v\right\|_{L^\ell}\le C \left\| e^{-tL}\sqrt{L}v\right\|_{L^q}
$$
provided 
$$
\frac1\ell+\frac{1}{2(n-1)}=\frac1q.
$$
Hence if $\ell<\frac{2n}{n-2}$, then $\cE_{0,\rho}(v)\in L^\ell(\bB(\rho))$ and
$$
\left\| \cE_{0,\rho}( v)\right\|_{L^\ell\left(\bB_0(\rho)\right)}^2\le
 C\rho^{\frac{2n}{\ell}}\, \left\| \sqrt{L}v\right\|_{L^2(S,\dvol_h)}^2
= C\rho^{\frac{2n}{\ell}-n+2}\,\langle\cN_{0,\rho} v,v\rangle_{L^2(\partial \bB_0(\rho), d \sigma_0)}\,\,.
$$

It is straightforward to deduce from all of this the more general result when $c_0\neq 0$: 
\begin {prop}
If $v\in H^{1/2}(\partial \bB_0(\rho))$ and $\ell<\frac{2n}{n-2}$, then
\begin{multline*}
\left\| \cE_{0,\rho} (v)\right\|_{L^\ell\left(\bB_0(\rho)\right)}^2\le  \\
C\rho^{\frac{2n}{\ell}-n+2}\,\langle\cN_{0,\rho} v,v\rangle_{L^2(\partial \bB_0(\rho), d \sigma_0)}+
C\rho^{\frac{2n}{\ell}-2n+2}\left(\int_{\partial \bB_0(\rho)} vd\sigma_0\right)^2.
\end{multline*}
\end {prop}
This estimate and the one in the next Proposition allow us to compare the spectra of $\cN_{0,\rho}$ and  $\cN_{V,\rho}$. 
\begin{prop}
 If $u\in W^{1,2}(\bB_0(\rho))$, then 
\begin{multline*}
\left(1-c\rho^{\bar\gamma}\right)\int_{\bB_0(\rho)} |du|^2_0 \dvol_0-C\rho^{\delta+1-n}
\left(\int_{\partial \bB_0(\rho)} u d\sigma_0\right)^2 \\ \le   \int_{\bB_0(\rho)}\left[\, |du|^2_g-Vu^2\,\right]\dvol_g \\
\le \left(1+c\rho^{\bar\gamma}\right)\int_{\bB_0(\rho)} |du|^2_0
\dvol_0+C\rho^{\delta+1-n}\left(\int_{\partial \bB_0(\rho)} u d\sigma_0\right)^2, 
\end {multline*}
where 
$\bar\gamma=\min\left\{\gamma, 2-\frac np\right\}$ and $\delta= 1-\frac{n}{p}$.\end{prop}
\proof
By hypothesis
$$
\left|\int_{\bB_0(\rho)}|du|^2_g \, \dvol_g-\int_{\bB_0(\rho)} |du|^2_0\right|\le 
C\Lambda\rho^{\gamma}\int_{\bB_0(\rho)} |du|^2_0 \, \dvol_0.
$$
Moreover, if $h=\cE_{0,\rho} R(u)$, then 
\begin{equation*}
 \begin{split}
  \left|\int_{\bB_0(\rho)}Vu^2\dvol_g\right|&\le 2\int_{\bB_0(\rho)}|V|(u-h)^2\dvol_g+2\int_{\bB_0(\rho)}|V|h^2\dvol_g\\
&\le C \|V\|_{L^{n/2}} \|u-h\|^2_{L^{\frac{2n}{n-2}}}+C\|V\|_{L^{p}}\|h\|^2_{L^{\frac{2p}{p-1}}}. 
 \end{split}
\end{equation*}
But $u-h\in W^{1,2}_0(\bB_0(\rho))$, so the Sobolev inequality and the variational characterization of $h$ yield
$$
\|u-h\|^2_{L^{\frac{2n}{n-2}}}\le C \|d(u-h)\|^2_{L^{2}}=C \left(\|du\|^2_{L^{2}}- \|dh\|^2_{L^{2}}\right)\le C\|du\|^2_{L^{2}}.
$$
Moreover, from the H\"older inequality,
$$
\|V\|_{L^{n/2}}\le C \rho^{2-\frac np}\|V\|_{L^{p}},
$$ 
and the previous proposition shows that 
$$
\|h\|^2_{L^{\frac{2p}{p-1}}}\le C\rho^{2-\frac n p} \langle\cN_{0,\rho} v,v\rangle_{L^2(\partial \bB_0(\rho), d \sigma_0)}+
C\rho^{1-\frac{n}{p}}\frac{1}{\rho^{n-1}}\left(\int_{\partial \bB_0(\rho)} u d\sigma_0\right)^2.
$$

Green's formula and the variational characterization of $h$ lead finally to
$$\langle\cN_{0,\rho} v,v\rangle_{L^2(\partial \bB_0(\rho), d \sigma_0)}=
\int_{\bB_0(\rho)} |dh|^2_0\dvol_0\le\int_{\bB_0(\rho)} |du|^2_0\dvol_0.$$
\endproof

One consequence of this proposition is that if $\rho$ is small enough, then $\cN_{V,\rho}$ has discrete spectrum
$$
\mu_0 <  \mu_1 \le\dots
$$
Moreover we obtain an estimate for the first two eigenvalues: 
\begin{prop}
$$
\left|\mu_0\right|\le C\rho^{1-\frac{n}{p}}, \qquad 
\left|\mu_1-\frac{\nu_1}{\rho}\right|\le C\rho^{\bar\gamma-1}. 
$$
\end{prop}

\subsection{A computation of eigenvalues}\label{spectral}
We now prove \eqref{nuSnuZ}. We begin with the identification $\displaystyle C(S)=\R^k\times C(Z) $ and recall
the form \eqref{sphersusp} of the metric $h$ on $S$.

First note that
\[
-\Delta_h =\bigoplus_{\mu \in \spec \Delta_Z \atop \lambda \in \spec (\Delta_{\bS^{k-1}})}L_{\mu,\lambda}  
\]
acting on $L^2\left( \left(0,\frac{\pi}{2}\right), \sin^{k-1} \psi\cos^{n-k-1}\psi\right)$, where 
\[
L_{\mu,\lambda}=-\frac{\partial^2\,}{\partial\psi^2}-\left((k-1)\cot \psi-(n-k-1)\tan\psi\right) 
\frac{\partial \, }{\partial\psi}+\frac{\mu}{\cos^2\psi}+\frac{\lambda}{\sin^2\psi}
\]
The first nonzero eigenvalue $-\Delta_h$ is the minimum of 
\begin{itemize}
\item the first non zero eigenvalue of $L_{0,0}$; 
\item the lowest eigenvalue of $L_{\mu_1,0}$, where $\mu_1$ is the first non zero eigenvalue of $-\Delta_Z$;
\item the lowest eigenvalue of $L_{0,\lambda_1}$, where $\lambda_1=k -1$ is the first non zero eigenvalue of $-\Delta_{\bS^{k-1}}$.
\end{itemize}

Now observe the following:
\begin{enumerate}
\item[i)] 
$\displaystyle L_{0,0}\left(\sin^2\psi-\frac{k}{n}\right)=2n \left(\sin^2\psi-\frac{k}{n}\right)$;
\item[ii)] Writing $\mu_1=\gamma(\gamma+n-k-2)$, then
$\displaystyle L_{\mu_1,0}\left(\cos^\gamma\psi\right)=\gamma(n-2+\gamma) \cos^\gamma\psi$; 
\item[iii)] $\displaystyle L_{0,k-1}\left(\sin\psi\right)=(n-1) \sin \psi$.
\end{enumerate}
These show that the first non zero eigenvalue of $-\Delta_h$ is 
\[
\begin{cases}
n-1& \,\,\mathrm{if}\,\, \mu_1 \ge n-k-1=\dim Z\\
\gamma(n-2+\gamma) & \,\,\mathrm{if}\,\, \mu_1 =\gamma(\gamma+n-k-2)\le n-k-1=\dim Z .
\end{cases}
\]

\section{Monotonicity formula}
Consider the truncated cone $C_R(S)$ with metric $g_0 = dr^2 + r^2 h$, where the link $S$ is a connected stratified space 
of dimension $n-1$ with an iterated edge metric $h$.  Consider another iterated edge metric $g$ which is Lipschitz 
with respect to $g_0$ and satisfies for all $\rho\in [r,R]$ :
$$
|g-g_0|\le \Lambda \rho^{\gamma} \ll 1 \quad \mathrm{on} \quad \bB_0(\rho)
$$

Since $S$ is connected, the spectrum of $-\Delta_h$ is a nondecreasing sequence 
$$\nu_0(n-2+\nu_0) = 0<\nu_1(n-2+\nu_1)\le\dots , $$
where 
$$
\nu_0=0<\nu_1\le \nu_2\le \dots \, .
$$
\begin{prop}\label{energy} Suppose that $V\in L^p$  for some $p > n/2$ and let $u\in W^{1,2}(C_R)$ satisfy
$$
\Delta_g u+Vu=0.
$$
Set $\bar\gamma=\min\left\{\gamma, 2-\frac n p \right\}$ and for any $\rho_-\le \rho_+\le R$ define
\[
\Psi(\rho_+, \rho_-) = 
\begin{cases}
\left|\rho_+^{2-n/p-2\nu_1}-\rho_-^{2-n/p-2\nu_1}\right| \quad & \mbox{if}\quad 1-\frac{n}{2p}-\nu_1\neq 0 \\
\log\left(\frac{\rho_+}{\rho_-}\right) \quad & \mbox{if} \quad 1-\frac{n}{2p}-\nu_1 =0.
\end{cases}
\]
Then there exists a constant $C$ depending only on $n,\Lambda,\nu_1, \|u\|_{L^\infty}$ and  $h$ such that
$$
\frac{e^{-C\rho_-^{\bar\gamma}}}{\rho_-^{n-2+2\nu_1}}\int_{\bB_0(\rho_-)}|du|^2_g\dvol_g\le
 \frac{e^{-C\rho_+^{\bar\gamma}}}{\rho_+^{n-2+2\nu_1}}\int_{\bB_0(\rho_+)}|du|^2_g\dvol_g+C\, \Psi(\rho_+,\rho_-).
$$
Moreover, there is a constant $\kappa$ such that if $\frac12 \tilde h\le h\le 2 \tilde h$ then 
$C(n,\nu_1, \Lambda,h)\le \kappa C(n,\nu_1, \Lambda,\tilde h)$.
\end{prop}
\proof We shall derive a differential inequality for the function
$$
\rho\mapsto E_0(\rho)= \int_{\bB_0(\rho)}|du|^2_0\dvol_0\,\,.
$$
First note that 
$$
E_0'(\rho)=\int_{\partial \bB_0(\rho)} |du|^2_0d\sigma_0=\int_{\partial\bB_0(\rho)} |d_Tu|^2_{\rho^2h}
d\sigma_0+\int_{\partial\bB_0(\rho)} \left| \frac{\partial u}{\partial r}\right|^2 d\sigma_0.
$$
where $d_T$ is the differential along $\partial \bB_0(\rho)$.  Next, 
$$
E(\rho)= \int_{\bB_0(\rho)}|du|^2_g\dvol_g=\int_{\bB_0(\rho)}Vu^2\dvol_g+\int_{\partial\bB_0(\rho)}u\cN_{V,\rho}ud\sigma_g
$$
satisfies 
$$
\left(1-c_n\Lambda \rho^\gamma\right)E( \rho)\le E_0(\rho)\le \left(1+c_n\Lambda \rho^\gamma\right)E( \rho).
$$

By \eqref{normal}, there is a constant $\eta$ such that 
\begin{multline*}
\left(1+\eta\Lambda\rho^\gamma\right)E_0'(\rho)-\frac{n-2+2\nu_1}{\rho }E(\rho) \\
\ge  \int_{\partial\bB_0(\rho)} |d_Tu|^2_{\rho^2h} d\sigma_0+\int_{\partial\bB_0(\rho)} \left| \cN_{V,\rho}u\right|^2 d\sigma_g\\
-\frac{n-2+2\nu_1}{\rho }\int_{\partial\bB_0(\rho)}u\cN_{V,\rho}ud\sigma_g
-\frac{n-2+2\nu_1}{\rho }\int_{\bB_0(\rho)}Vu^2\dvol_g\,.
\end{multline*}

We now compare $\int_{\partial\bB_0(\rho)}u\cN_{V,\rho}ud\sigma_g$ and $\int_{\partial\bB_0(\rho)}u\cN_{0,\rho}ud\sigma_0$.
Introducing the harmonic  $h:=\cE_{0,\rho}\left(\left.u\right|_{\partial \bB_{0}(\rho)}\right)$, then we have
\begin{equation*}\begin{split}
\int_{\partial\bB_0(\rho)}u\cN_{V,\rho}ud\sigma_g&-\int_{\partial\bB_0(\rho)}u\cN_{0,\rho}ud\sigma_0\\
&=\int_{\bB_0(\rho)}|du|^2_g\dvol_g-\int_{\bB_0(\rho)}Vu^2\dvol_g- \int_{\bB_0(\rho)}|dh|^2_0\dvol_0\,.
\end{split}\end{equation*}
Since $u\in L^\infty$, 
$$
\left|\int_{\bB_0(\rho)}Vu^2\dvol_g\right|\le C \int_{\bB_0(\rho)}|V|\le C \rho^{n\left(1-\frac 1p\right)},
$$
and moreover, by the variational characterization of $h$, 
$$
\int_{\bB_0(\rho)}|dh|^2_0\dvol_0\le \int_{\bB_0(\rho)}|du|^2_0\, \dvol_0\le\left(1+c\Lambda\rho^\gamma\right) \int_{\bB_0(\rho)}|du|^2_g\, \dvol_g.
$$
Using the same argument for $u$ and the fact that $\|h\|_{L^\infty}\le \|u\|_{L^\infty}$, we get
\begin{multline*}
\int_{\bB_0(\rho)}|du|^2_g\, \dvol_g-\int_{\bB_0(\rho)}Vu^2 \, \dvol_g \le \int_{\bB_0(\rho)}|dh|^2_g\dvol_g-
\int_{\bB_0(\rho)}Vh^2\, \dvol_g\\
\le  \left(1+c\Lambda\rho^\gamma\right)\int_{\bB_0(\rho)}|dh|^2_0\, \dvol_0+ C \rho^{n\left(1-\frac 1p\right)}.
\end{multline*}
Hence there is a constant depending only on $V$, $n$ and $\|u\|_{L^\infty}$ such that 
\begin{equation}
\left|\int_{\partial\bB_0(\rho)}u\cN_{V,\rho}u\, d\sigma_g-\int_{\partial\bB_0(\rho)}u\cN_{0,\rho}u \, d\sigma_0\right|\le 
C\Lambda \rho^\gamma E_0(\rho)+C\rho^{n\left(1-\frac 1p\right)}.
\end{equation}
But 
\begin{multline*}
\int_{\partial\bB_0(\rho)} |d_Tu|^2_{\rho^2h}
d\sigma_0 =\frac{1}{\rho^2} \int_{\partial\bB_0(\rho)} u\Delta_hud\sigma_0\\
= \int_{\partial\bB_0(\rho)} \cN_{0,\rho}u\left(\cN_{0,\rho}u+\frac{n-2}{\rho}u\right)d\sigma_0\\
\ge \frac{n-2+\nu_1}{\rho} \int_{\partial\bB_0(\rho)}u \cN_{0,\rho}u\,\,d\sigma_0
\end{multline*}
so that 
\begin{equation*}\begin{split}
\int_{\partial\bB_0(\rho)} |d_Tu|^2_{\rho^2h}d\sigma_0&+\int_{\partial\bB_0(\rho)} \left| \cN_{V,\rho}u\right|^2 d\sigma_g-\frac{n-2+2\nu_1}{\rho }\int_{\partial\bB_0(\rho)}u\cN_{V,\rho}ud\sigma_g\\
&\ge \frac{n-2+\nu_1}{\rho}\left[ \int_{\partial\bB_0(\rho)}u \cN_{0,\rho}u\,\,d\sigma_0
- \int_{\partial\bB_0(\rho)}u \cN_{V,\rho}u\,\,d\sigma_g\right]\\
&+\int_{\partial\bB_0(\rho)} \left| \cN_{V,\rho}u\right|^2 d\sigma_g-\mu_1   \int_{\partial\bB_0(\rho)}u \cN_{V,\rho}u\,\,d\sigma_g   \\
&+\left(\mu_1-\frac{\nu_1}{\rho}\right)  \int_{\partial\bB_0(\rho)}u \cN_{V,\rho}u\,\,d\sigma_g.
\end{split}\end{equation*}

By our comparison result above, the first term on the left is bounded from below by
$$-C\left(\rho^{\gamma-1}\right)E_0(\rho)-C\rho^{n\left(1-\frac1p\right)-1}.$$
Using the spectral theorem and the estimate on the eigenvalues of the Dirichlet to Neumann operator $\cN_{V,\rho}$,  the second term in the LHS is bounded from below by 
$$|\mu_0\mu_1|\,\int_{\partial\bB_0(\rho)} u^2\ge -C \rho^{n-2+\delta}=C\rho^{n\left(1-\frac1p\right)-1}$$
Similarly, the last term in the LHS is bounded from below by
$$-C\left(\rho^{\bar\gamma-1}\right)E_0(\rho)-C\rho^{n\left(1-\frac1p\right)-1}.$$
Eventually, we get a constant $\kappa$ such that we have the differential inequality
$$\left(1+\kappa\rho^{\bar\gamma}\right)E_0'(\rho)-\frac{n-2+2\nu_1}{\rho }E_0(\rho)\ge -C\rho^{n\left(1-\frac 1p\right)-1}. 
$$
The result follows now easily. \endproof

\section{Proof of \tref{A}}\label{proofs}
We now turn to a proof of our first main theorem. Let $M$ be an $n$-dimensional stratified space with an iterated edge metric $g$. 
Assume that each unit tangent sphere $S_m$, $m\in M$ is connected and that for some $\nu\in (0,1],$ for all $m\in M$ the 
first nonzero eigenvalue of the Laplace operator on $S_m$ is larger than $\nu(n-2+\nu)$. Suppose that $V\in L^p$ for some $p>n/2$ and  
$u\in W^{1,2}(M)$ a solution of the equation $\displaystyle \Delta u + Vu = 0$. We know already that $u\in L^\infty$, and our
goal is to show that $u$ has a certain H\"older regularity. 

\medskip
\noindent {\bf First case: $\nu=1$ and $V\in L^\infty$: }
By \pref{energy} and \tref{geoBoule}, we see that for all  $p\in M$ and $r\in (0,\eta)$,
$$
\frac{e^{-Cr^{\bar\gamma}}}{r^{n}}\int_{\bB_0(r)}|du|^2_g\le C+C|\log(r)|.
$$
The second remark in the appendix shows that there is a constant $C$ such that for all $x,y\in M$:
$$
\left|u(x)-u(y)\right|\le C\sqrt{|\log(d(x,y))|} d(x,y)\,.
$$

\medskip

\noindent {\bf Second case: $\nu<1$ and $V\in L^\infty$:}
According to \pref{energy} and \tref{geoBoule}, if $p\in M$ and $r\in (0,\eta)$, then 
$$
\frac{e^{-Cr^{\bar\gamma}}}{r^{n-2+2\nu}}\int_{\bB_0(r)}|du|^2_g\le C.
$$
Applying the H\"older result Proposition~\ref{morrey} in this setting proves the result.

\medskip

\noindent {\bf Third case: $V\in L^p$ where $p\in (n/2,\infty)$:}
In this case, by \pref{energy} and \tref{geoBoule}, we obtain that for all $p\in M$ and $r\in (0,\eta)$ :
$$
\frac{e^{-Cr^{\bar\gamma}}}{r^{n-2+2\nu}}\int_{\bB_0(r)}|du|^2_g\le C+Cr^{2-\frac{n}{p}-2\nu}
$$
so that 
$$
\int_{\bB_0(r)}|du|^2_g\le Cr^{n-2+2\nu}+Cr^{n-\frac{n}{p}}.
$$
Hence if we set $\mu=\min\left\{\nu, 1-\frac{n}{2p}\right\}$ then $\displaystyle \int_{\bB_0(r)}|du|^2_g\le Cr^{n-2+2\mu}$, so
by Proposition~\ref{morrey} again, $u$ is H\"older continuous of order $\mu$. 

\medskip

\noindent {\bf Some remarks: } 
\begin{enumerate}[i)]
\item Suppose that $\tilde g$ is another iterated edge metric on $M$ such that $\tilde g-g=\sigma^{\gamma}h$, where $h$ is 
an iterated edge symmetric two tensor, $\sigma$ is the distance to $M^{\mathrm{sing}}$ and  $\gamma>0$. Then solutions 
of the equation $\Delta_{\tilde g}u+Vu=0$ have the same H\"older regularity as for the corresponding equation relative
to the metric $g$. 
\item We have seen that if $u\in W^{1,2}$ satisfies $\Delta u\in L^p$ for some $p>n/2$, then $u\in L^{\infty}$.  
So if we define $v=u+2\|u\|_{L^\infty}$, then $v$ is a solution of 
$$
\Delta v+Vv=0, 
$$ 
where 
$$
V=-\frac{\Delta u}{v}\in L^p.
$$
This means that $v$, and hence $u$, are also H\"older continuous of order $\mu$. 
\item A point has capacity zero, so the equation $\Delta u\in L^p$ also holds when we remove a finite number of points from
$M$. By  \rref{cut}, if the condition on  the connectedness of the spheres is not satisfied then $u$ is H\"older continuous 
of order $\alpha$ on $\overline{M}$.

\item There are general results for the regularity of solution of the equation $\Delta u\in L^p$ on Metric Measure space;
for instance \cite{Jiang,Ko} contains a result about Lipschitz continuity of solutions under the condition that the underlying
measure is Alfors regular, and that there is a uniform Poincar\'e inequality and a kind of heat kernel-curvature lower bound. 
In our setting, even harmonic function may not Lipschitz, and our result are optimal with respect to the exponent of regularity.
\end{enumerate}

\appendix
\section{Appendix: Morrey implies H\"older}
In this appendix, we recall the proof that a Morrey-type regularity  H\"older regularity result based on Morrey's idea. 

Suppose that $(M,d,\mu)$ is a compact almost smooth metric-measure space which satisfies the following properties: 
\begin{enumerate}
\item[i)]  $d\mu$ is a doubling measure, i.e.\ there is some $V>0$ such that 
\[
\mu (B(p,2r)) \le V \mu (B(p,r)) 
\]
for every point $p \in M$ and $r < \mbox{diam}(M)/2$. 
\item[ii)] The uniform Poincar\'e inequality holds: there exist $ A\ge 1 $ and $C > 0$  such that
\begin{equation}
\|f-f_B\|^2_{L^2(B(p,r))} \le C r^2 \int_{B(p,Ar) }|df|^2 \, d\mu
\label{Poincare}
\end{equation}
for all $f \in W^{1,2}(B(p,Ar); d\mu)$,  $p \in M$ and $r < \mbox{diam}(M)/(2A)$. 
(Here $f_B:=\frac{1}{\mu (B(p,r))}\int_{B(p,r)} f\, d\mu$.) 
\end{enumerate}

\begin{prop}\label{morrey} Assume that $v\in W^{1,2}(M; d\mu)$ satisfies 
\[
\frac{1}{\mu (B(p,r))}\int_{B(p,r) }|dv|^2 \, d\mu\le \Lambda r^{2\alpha-2}
\] 
for some $\Lambda>0$, $\alpha\in (0,1]$ and $\eta>0$, and for every $p\in M$ and $r\in (0,\eta) $. 
Then $v$ is $\alpha$-H\"older continuous. (In the special case $\alpha  = 1$, we mean that $v$ is Lipschitz.)
\end{prop}
\proof The proof is classical, see for instance \cite[Lemme 3.4]{CarrobJga} for other applications of these ideas. 

First note that if $\ell\in \N$ is chosen so that $A\le 2^\ell$, then 
\[
\mu (B(p,2Ar)) \le V^{\ell+1} \mu (B(p,r)\,\,. 
\]
For $p\in M$, suppose that $2Ar<\eta$, and write $B:=B(p,r)$, $2B:=B(p,2r)$, and 
$$
v_B=\frac{1}{\mu (B)}\int_{B} vd\mu, \,\, v_{2B}=\frac{1}{\mu (2B)}\int_{2B} v\, d\mu  \,\,.
$$ 
Then
 \begin{equation*}
 \begin{split}
 \left|v_B-v_{2B}\right|&=\frac{1}{\mu (B)\mu (2B)}\left|\int_{B\times2B} 
 \left( v(x)-v(y) \right)d\mu(x)d\mu(y) \right|\\
 &\le \frac{1}{\sqrt{\mu (B)\mu (2B)}}\left(\int_{B\times2B} 
 \left( v(x)-v(y) \right)^2d\mu(x)d\mu(y) \right)^{\frac12}\\
 &\le  \frac{1}{\sqrt{\mu (B)\mu (2B)}}\left(\int_{2B\times2B} 
 \left( v(x)-v(y) \right)^2d\mu(x)d\mu(y) \right)^{\frac12}\\
 &\le  \frac{\sqrt{2}}{\sqrt{\mu (B)}}\left(\int_{2B} 
 \left( v-v_{2B}\right)^2d\mu\right)^{\frac12}\\
 &\le  \frac{\sqrt{2C}}{\sqrt{\mu (B)}}2r \left(\int_{B(p,2Ar)} 
 |dv|^2d\mu \right)^{\frac12}\\
 &\le \frac{\sqrt{2C\Lambda}}{\sqrt{\mu (B)}}2r \sqrt{\mu (B(p,2Ar))}(2A)^{\alpha-1}r^{\alpha-1}\\
 &\le \sqrt{8V^{\ell+1}C\Lambda }\, (2A)^{\alpha-1} r^\alpha.
 \end{split}\end{equation*}
 
For any $\rho\in (0,\eta/(2A))$, apply this inequality to $r=\rho/2^{k}$, $k=1,2,\dots$ and sum the inequalities over all $k$.
We obtain in this way a constant $\kappa>0$ such that for any $\rho\in (0,\eta/(2A))$ and any $p\in M$,
\[
\left|v(p)-v_{B(p,\rho)}\right|\le \kappa \Lambda \rho^\alpha
\]
Hence if $4Ad(x,y)\le \eta$, then 
\begin{multline*}
\left|v(x)-v(y)\right|\\
\le \left|v(x)-v_{B(x,d(x,y))}\right|+\left|v(y)-v_{B(y,d(x,y))}\right|+\left|v_{B(x,d(x,y))}-v_{B(y,d(x,y))}\right|\\
\le 2\kappa d(x,y)^\alpha+\left|v_{B(x,d(x,y))}-v_{B(y,d(x,y))} \right|. 
\end{multline*}

The same argument gives that for $d:=d(x,y)$, 
\begin{equation*}
 \begin{split} &\left|v_{B(x,d)}-v_{B(y,d)} \right|\\
   &\le\frac{1}{\sqrt{\mu (B(x,d))\mu (B(y,d))}}\left(\int_{B(x,2d)\times B(x,2d) } 
 \left( v(t)-v(z) \right)^2d\mu(t)d\mu(z) \right)^{\frac12}\\
   &\le\frac{\sqrt{2\mu (B(x,2d))}}{\sqrt{\mu (B(x,d))\mu (B(y,d))}}\left(\int_{B(x,2d) } 
 \left( v-v_{B(x,2d)} \right)^2d\mu \right)^{\frac12}\\
 &\le \kappa' \Lambda d^\alpha \sqrt{\frac{\mu (B(x,2Ad))}{\mu (B(y,d))}}\\
 &\le \kappa' \Lambda d^\alpha \sqrt{\frac{\mu (B(y,(2A+1)d))}{\mu (B(y,d))}}\\
 &\le \kappa' \Lambda d^\alpha V^{\ell+2}.
\end{split}\end{equation*}
This proves the result.
\endproof

\begin{rems} \label{remarapp}
This argument is local. Hence if $v\in W^{1,2}(\Omega; d\mu)$ satisfies 
\[
\frac{1}{\mu (B)}\int_{B }|dv|^2 \, d\mu\le \Lambda r(B)^{2-2\alpha}
 \]
for all balls $B\subset \Omega$ of radius $r(B)\in (0,\eta)$, then 
$v$ is $\alpha$-H\"older continuous on $\Omega$. In fact, for any $\delta>0$ 
there is a constant $\kappa$ such that if 
$x,y\in \Omega$ and $d(x,\partial \Omega)\ge \delta$, $d(y,\partial \Omega)\ge \delta$, then
\[
\left|v(x)-v(y)\right|\le \kappa d(x,y)^\alpha.
\]

It is also easy to check that if $u$ satisfies 
\[
\frac{1}{\mu (B(p,r))}\int_{B(p,r) }|dv|^2 \, d\mu\le \Lambda \left|\log(r)\right|^{2\gamma},
\]
for all $p\in M$, $r\in (0,\eta)$ and for some $\gamma>0$, then there is a constant $C$ such that 
\[
|u(x)-u(y)|\le C \left|\log(d(x,y)\right|^{\gamma}d(x,y) \quad \forall\, x,y \in M.
\]
\end{rems}

\end{document}